\documentclass[11pt]{article}

\usepackage{amscd,amsmath, amssymb, wrapfig, epsfig}

\newcommand{\bftext}[1]{{\bf #1}}

\newcommand{\RR}{\mathbb{R}}

\numberwithin{equation}{section}

\def\eqref#1{(\ref{#1})}
\newcommand{\goth}{\mathfrak}

\newcommand{\arrow}{{\:\longrightarrow\:}}

\newcommand{\C}{{\Bbb C}}
\newcommand{\R}{{\Bbb R}}

\newcommand{\6}{\partial}
\def\1{\sqrt{-1}\:}
\newcommand{\restrict}[1]{{\left|_{{\phantom{|}\!\!}_{#1}}\right.}}

\newcommand{\calo}{{\cal O}}
\newcommand{\cac}{{\cal C}}

\renewcommand{\tilde}{\widetilde}
\renewcommand{\bar}{\overline}
\renewcommand{\phi}{\varphi}
\renewcommand{\epsilon}{\varepsilon}
\renewcommand{\geq}{\geqslant}
\renewcommand{\leq}{\leqslant}

\newcommand{\Lie}{\operatorname{Lie}}
\newcommand{\End}{\operatorname{End}}
\newcommand{\Tot}{\operatorname{Tot}}
\newcommand{\Id}{\operatorname{Id}}

\newcommand{\Aut}{\operatorname{Aut}}

\newcommand{\Spec}{\operatorname{Spec}}

\renewcommand{\Re}{\operatorname{Re}}

\newcounter{Mycounter}[section]
\newcounter{lemma}[section]
\newcounter{claim}[section]
\renewcommand{\theclaim}{{Claim \thesection.\arabic{claim}}}
\newcommand{\claim}{%
     \setcounter{claim}{\value{Mycounter}}
     \refstepcounter{claim}
     \stepcounter{Mycounter}
     {\noindent \bf \theclaim.\ }}

\newcounter{sublemma}[section]
\newcounter{corollary}[section]
\renewcommand{\thecorollary}{{Corollary
\thesection.\arabic{corollary}}}
\newcommand{\corollary}{%
     \setcounter{corollary}{\value{Mycounter}}
     \refstepcounter{corollary}
     \stepcounter{Mycounter}
     {\noindent \bf \thecorollary.\ }}

\newcounter{theorem}[section]
\renewcommand{\thetheorem}{{Theorem \thesection.\arabic{theorem}}}
\newcommand{\theorem}{%
     \setcounter{theorem}{\value{Mycounter}}
     \refstepcounter{theorem}
     \stepcounter{Mycounter}
     {\noindent \bf \thetheorem.\ }}

\newcounter{conjecture}[section]
\newcounter{proposition}[section]
\renewcommand{\theproposition}
       {{Proposition \thesection.\arabic{proposition}}}
\newcommand{\proposition}{%
     \setcounter{proposition}{\value{Mycounter}}
     \refstepcounter{proposition}
     \stepcounter{Mycounter}
     {\noindent \bf \theproposition.\ }}

\newcounter{definition}[section]
\renewcommand{\thedefinition}
       {{Definition~\thesection.\arabic{definition}}}
\newcommand{\definition}{%
     \setcounter{definition}{\value{Mycounter}}
     \refstepcounter{definition}
     \stepcounter{Mycounter}
     {\noindent \bf \thedefinition.\ }}

\newcounter{example}[section]
\renewcommand{\theexample}{{Example \thesection.\arabic{example}}}
\newcommand{\example}{%
     \setcounter{example}{\value{Mycounter}}
     \refstepcounter{example}
     \stepcounter{Mycounter}
     {\noindent \bf \theexample:\ }}

\newcounter{remark}[section]
\renewcommand{\theremark}{{Remark \thesection.\arabic{remark}}}
\newcommand{\remark}{%
     \setcounter{remark}{\value{Mycounter}}
     \refstepcounter{remark}
     \stepcounter{Mycounter}
     {\noindent \bf \theremark.\ }}

\newcounter{problem}[section]
\newcounter{question}[section]
\renewcommand{\thequestion}{{Question
\thesection.\arabic{question}}}
\newcommand{\question}{%
     \setcounter{question}{\value{Mycounter}}
     \refstepcounter{question}
     \stepcounter{Mycounter}
     {\noindent \bf \thequestion.\ }}

\makeatletter

\newcommand{\ps@verbit}{%
  \renewcommand{\@oddhead}{%
          \scriptsize
          {Sasakian structures on CR-manifolds}
          \hfil\tiny {L. Ornea and M. Verbitsky, June 6, 2006}}
  \renewcommand{\@evenhead}{\@oddhead}
  \renewcommand{\@oddfoot}{\hfil\thepage\hfil}
  \renewcommand{\@evenfoot}{\@oddfoot}}

\@addtoreset{equation}{section}
\@addtoreset{footnote}{section}
\makeatother

\def\blacksquare{\hbox{\vrule width 5pt height 5pt depth 0pt}}
\def\endproof{\blacksquare}

\begin{document}
\begin{center}
{\LARGE\bf
Sasakian structures on CR-manifolds\\[3mm]
}

Liviu Ornea and Misha Verbitsky\footnote{Misha Verbitsky is an EPSRC
advanced
fellow supported by  EPSRC grant
GR/R77773/01.

{ {\bf Keywords:} CR-manifold, Sasakian manifold, Reeb field, pseudo-convex
manifold, Vaisman
manifold, Stein manifold, deformation, potential.}

\scriptsize
{\bf 2000 Mathematics Subject
Classification:} { 53C55}}

\end{center}

{\small
\hspace{0.15\linewidth}
\begin{minipage}[t]{0.7\linewidth}
{\bf Abstract} \\
A contact manifold $M$ can be defined as a quotient of a
symplectic manifold $X$ by a proper, free action of $\R^{>0}$,
with the symplectic form homogeneous of degree 2. If $X$
is, in addition, K\"ahler, and its metric is also homogeneous 
of degree 2, $M$ is called Sasakian. A Sasakian manifold is
realized naturally as a level set of a K\"ahler potential
on a complex manifold, hence it is equipped with a
pseudoconvex CR-structure. We show that any Sasakian
manifold $M$ is CR-diffeomorphic to an $S^1$-bundle of
unit vectors in a positive line bundle on a projective
K\"ahler orbifold. This induces an embedding from $M$
to an algebraic cone $C$. We show that this embedding
is uniquely defined by the CR-structure. Additionally, we
classify the Sasakian metrics on an odd-dimensional sphere
equipped with a standard CR-structure.
\end{minipage}
}

\tableofcontents


\section{Introduction}


\subsection{Sasakian manifolds and algebraic cones}

In this paper, we study existence of
Sasakian metrics on strictly pseudoconvex
CR-manifolds.
A pseudoconvex CR-manifold is a geometric
structure arising on a smooth boundary of a Stein
domain $X$ (\ref{_CR_stri_pse_Definition_},
\ref{_Stein_CR_Remark_}). If $M$ is compact and strictly pseudoconvex
CR-manifold of dimension $>3$,
then $M$ can always be realized as
a boundary of a Stein variety $X$ with at most
isolated singularities (\cite{andreotti_siu},
\cite{_Marinesku_Yeganefar_}). In fact, the
geometry of CR-structures is essentially
the same as the holomorphic geometry of the
corresponding Stein variety. In particular,
the automorphisms of $X$ are in a natural
correspondence with the CR-diffeomorphisms
of its boundary.

Strictly pseudoconvex CR-manifolds are always contact;
they are sometimes called {\em contact pseudoconvex}.

Sasakian metrics are the special  Riemannian
metrics on contact  pseudoconvex CR-manifolds.
They are related to K\"ahler metrics, in the same
way as the contact structures are related to symplectic
structures. A contact manifold can be defined
as a manifold with a symplectic structure
on its cone; a Sasakian metric on a contact
manifold induces a K\"ahler metric on
its symplectic cone (\ref{_Sasakian_Definition_}).

Sasakian manifolds can be defined in terms of algebraic
cone spaces, as follows.

\hfill

\definition\label{_algebra_cone_intro_Definition_}
{\bf A closed algebraic cone}  is an
affine variety $\cac$ admitting a $\C^*$-action $\rho$
with a unique fixed point $x_0$, which satisfies
the following.

\begin{description}
\item[1.] $\cac$ is smooth outside of $x_0$.
\item[2.] $\rho$ acts on the Zariski tangent
space $T_{x_0}\cac$ diagonally, with all eigenvalues
$|\alpha_i|<1$.
\end{description}

\noindent {\bf An open algebraic cone} is $\cac \backslash \{x_0\}$.

\hfill

In Section \ref{_cones_Section_} we give another, equivalent
but more constructive, definition of an algebraic cone
(\ref{_alge_cone_clo_ope_Definition_}).

\hfill

By definition, a Sasakian manifold $M$ admits
a CR-embedding into an algebraic cone $\cac(M):= M\times \R^{>0}$,
as a set $M\times \{t_0\}$. The function
$\cac(M) \arrow \R^{>0}$, $(m, t)\arrow t^2$ is
a K\"ahler potential of $\cac(M)$, as follows
from an elementary calculation (see \emph{e.g.}
\cite{_Verbitsky:LCHK_}). The converse is also true:
given a K\"ahler potential $\phi:\; \cac \arrow \R$
on an algebraic cone $\cac$, satisfying
$\Lie_v\phi = 2\phi$, for a vector field $v\in T\cac$
inducing a holomorphic contraction on $\cac$,
we may assume that $(\cac, \6\bar\6\phi)$ is
a Riemannian cone of $M$.\footnote{Here, $\Lie_v$ denotes the Lie derivative.}

The correspondence between algebraic cones and
Sasakian manifolds is quite significant. One may argue
that the algebraic cone, associated
with a Sasakian manifold, gives a functor
similar in many respects to the forgetful
functor from the category of K\"ahler manifolds
to the category of complex manifolds.
Indeed, the moduli space of algebraic cones
is finite-dimensional, as follows from
\ref{_alge_cone_clo_ope_Definition_}, and the Sasakian metrics
are determined by an additional set of $C^\infty$-data
(the K\"ahler potential).

One could also argue that a proper analogy
of a complex structure is a CR-structure underlying a Sasakian
manifold. However, the CR-structure (unlike
complex structure, or a structure of an algebraic cone)
in many cases, \emph{e.g.} in dimension 3, determines the
Sasakian metric completely (up to a constant).
 In fact, there is only a finite-dimensional set of
Sasakian metrics on a given CR-manifold (\cite{bgs};
see also \ref{_when_CR_admits_Sasakian_Theorem_}).

In this paper, we study the forgetful functor
from the category of Sa\-sa\-kian manifolds
to the category of algebraic cones. We show
that it is determined by the CR-structure.

\hfill

\theorem\label{_S^1_equiv_Sasakian_Theorem_}
Let $M$ be a compact pseudoconvex contact CR-manifold.
Then the following conditions are equivalent.
\begin{description}
\item[(i)] $M$ admits a Sasakian metric, compatible
with the CR-structure.
\item[(ii)] $M$ admits a proper, transversal
CR-holomorphic $S^1$-action.
\item[(iii)] $M$ admits a nowhere degenerate, transversal
CR-holomorphic vector field.
\end{description}

\hfill

\theorem\label{_Sasa_unique_Theorem_}
Let $M$ be a compact, strictly pseudoconvex CR-manifold
admitting a proper, transversal CR-\-ho\-lo\-mor\-phic $S^1$-\-ac\-tion.
Then $M$ admits a unique (up to an automorphism)
$S^1$-invariant CR-embedding into an algebraic cone $\cac$.
Moreover, a Sasakian metric on $M$ can be induced
from an automorphic K\"ahler metric on this cone.

\hfill

We prove \ref{_S^1_equiv_Sasakian_Theorem_}
in Subsection \ref{exist_subsection},
and \ref{_Sasa_unique_Theorem_} in Subsection
\ref{uni_subsection}  (when $M$ is not a sphere).
The case when $M$ is a sphere is considered at
the end of Subsection \ref{_Sasa_on_sphe_Subsection_}.

\hfill

\remark
The Sasakian metric is by definition induced from
an embedding to its cone, which is a K\"ahler manifold.
This cone is algebraic, as indicated above.
This metric is not unique, though the embedding
is unique and canonical, as follows from
\ref{_Sasa_unique_Theorem_}.

\hfill

\remark \label{_another_arguments_citation_Remark_}
For another approach to the existence of Sasakian structures
compatible with a contact pseudoconvex structure on a compact
manifold, see \cite{bgs}.
A still different approach is the following: In the course of the proof of
\cite[Theorem E]{lee}, Lee proves that the infinitesimal generator
of a transverse CR-automorphism (of a pseudoconvex contact
structure) is necessarily a Reeb vector field for a contact form
underlying the given contact bundle. On the other hand, Webster
proved in \cite{webster} that if the Reeb field of a pseudoconvex
contact structure is a CR-automorphism, then the torsion of the
Tanaka connection vanishes. But it is known (see \emph{e.g.}
\cite{sorin}) that a pseudoconvex contact structure with zero Tanaka
torsion is Sasakian. However, it seems that this result was
never explicitely stated as such.

\hfill

\remark
In dimension 3, Sasakian structures on CR-manifolds were
completely classified (\cite{geiges}, \cite{belgun}, \cite{_Belgun:S^3_}).
For a 3-dimensional CR-manifold $M$, not isomorphic to a sphere,
the Sasakian metric  is unique, hence the corresponding
cone is also unique. \ref{_S^1_equiv_Sasakian_Theorem_}
and \ref{_Sasa_unique_Theorem_} in dimension 3
follow immediately from \cite{belgun}, \cite{_Belgun:S^3_}.
In this paper, we shall always assume that
$\dim M \geq 5$.

\subsection{Sasakian geometry and contact geometry}

There is a way to define contact manifolds and Sasakian
manifolds in a uniform manner. Let $M$ be a smooth
manifold equipped with a free, proper  action
of the multiplicative group $\R^{>0}$, and
a symplectic form $\omega$. Assume that
$\omega$ is homogeneous of weight 2 with respect to
$\rho$, that is, $\Lie_v\omega= 2\omega$,
where $v\in TM$ is the tangent vector field
of $\rho$. Then the quotient
$M/\rho$ is contact. This can be considered
as a definition of a contact manifold
(see \ref{_conical_symple_Remark_}).
Then $M$ is called a \bftext{symplectic cone}
of a contact manifold $M/\rho$.

Now, let  $(M, g, \omega)$ be a K\"ahler manifold
(here we consider $(M,\omega)$ as a symplectic manifold,
equipped with a compatible Riemannian structure $g$).
Assume again that $\rho$ is a free, proper action
of $\R^{>0}$ on $M$, and $g$ and $\omega$ are homogeneous
of weight 2:
\[
\Lie_v\omega= 2\omega, \ \  \Lie_v g= 2g.
\]
The quotient $M/\rho$ is contact (as indicated above)
and Riemannian (the Riemannian metric is obtained
from $g$ by appropriate rescaling). It easily
follows from the definitions that $M/\rho$ is
Sasakian. In fact, the Sasakian manifolds can
be defined this way (see \ref{_Sasakian_Definition_}).
The Sasakian metric
is therefore a natural odd-dimensional
counterpart to K\"ahler metrics.

In symplectic geometry, one is often asked the
following question.

\hfill

\question
Let $M$ be a symplectic manifold. Is there a K\"ahler
metric compatible with the symplectic structure?

\hfill

It is natural to ask the same question for contact
geometry.

\hfill

\question
Let $M$ be a contact manifold. Is there a Sasakian
metric compatible with the contact structure?

\hfill

A partial answer to this question is given in this
paper, in the additional assumptions of existence of CR-structure,
which is natural and very common
in contact topology.

\hfill

A set of natural examples of CR and Sasakian
manifolds is provided by algebraic geometry.

\hfill

\example\label{_U(1)_fibration_Examples_}
Let $X$ be a projective orbifold with quotient
singularities (in algebraic geometry such an
object is also known under the name of ``Deligne-Mumford stack''),
and $L$ an ample Hermitian line bundle on $X$. Assume
that the curvature of $L$ is positive.
Let $\Tot(L^*)$ be the space of all non-zero vectors
in the dual bundle, considered as a complex manifold,
and $\phi:\; \Tot(L^*)\arrow \R$
map $v\in L^*$ into $|v|^2$. It is easy to check
that $\phi$ is strictly plurisbuharmonic, that is,
$\6\bar\6\phi$ is a K\"ahler form on $\Tot(L^*)$.
Therefore, the level set $M:= \phi^{-1}(\lambda)$
of $\phi$ is a  strictly pseudoconvex CR-manifold.
This level set is a $U(1)$-bundle on $X$.

It is easy to see that the metric $\6\bar\6\phi$
induces a Sasakian structure on $M$
(see e.g. \cite{_Verbitsky:LCHK_}).
Such Sasakian manifolds are called \bftext{quasiregular}.

\hfill

In dimension 3, F. A. Belgun has shown that all
Sasakian manifolds are obtained this way (see \cite{geiges}, \cite{belgun},
\cite{_Belgun:S^3_}):

\hfill

\theorem
A strictly pseudoconvex, compact
 CR-manifold $M$ of dimension
3 admits a Sasakian metric if and only if $M$
is isomorphic to a $U(1)$-fibration associated with a
positive line bundle on a projective orbifold
(\ref{_U(1)_fibration_Examples_}).
Moreover, the Sasakian metric on $M$
is unique, up to a constant multiplier, unless
$M$ is $S^3\subset \C^2$.

\hfill

Using a construction of Sasakian positive
cone due to \cite{bgs}, we shall extend this theorem
to arbitrary dimension.

\hfill

\theorem\label{_when_CR_admits_Sasakian_Theorem_}
Let $M$ be a strictly pseudoconvex, compact
 CR-manifold. Then $M$ admits a Sasakian metric
if and only if $M$
is CR-isomorphic to a $U(1)$-fibration associated with a
positive line bundle on a projective orbifold
(\ref{_U(1)_fibration_Examples_}).
Moreover, the set of Sasakian structures
on $M$ is in bijective correspondence with the
set of positive and transversal
CR-holomorphic vector fields on
$M$.\footnote{This set is called {\bf the positive
Sasakian cone of $M$}. For a definition of
positive and transversal CR-holomorphic
vector fields, see Subsection \ref{_posi_Sasa_cone_Subsection_}.}

\hfill

\noindent{\bf Proof:} The last claim of
\ref{_when_CR_admits_Sasakian_Theorem_}
is due to \cite{bgs}; we give a new proof
of this statement in Subsection
\ref{_posi_Sasa_cone_Subsection_},
and prove the rest of
\ref{_when_CR_admits_Sasakian_Theorem_}.
\endproof


\section{Strictly pseudoconvex CR-manifolds}


\subsection{CR-manifolds and contact manifolds}

We recall some definitions, which are well known.

\hfill

\definition
Let $M$ be a smooth manifold. A \bftext{ CR-structure}
(Cauchy-Riemann structure) on $M$ is a subbundle
$H\subset TM\otimes \C$ of the complexified tangent bundle,
which is closed under commutator:
\[
[H,H]\subset H
\]
and satisfies $H\cap \bar H=0$.

\hfill

A complex manifold $(P,I)$ is considered as a CR-manifold,
with $H=T^{1,0}P\subset TP\otimes \C$.

\hfill

\definition
Consider a real submanifold $M\subset P$ in a complex
manifold $(P,I)$. Suppose that $TM\cap I(TM)$ has constant
rank. Clearly, $H_M:=TM\otimes \C \cap T^{1,0}P$
is a CR-structure on $M$. Then $(M, H_M)$
is called \bftext{a CR-submanifold in $M$},
and $H_M$ is called \bftext{the induced CR-structure}.

\hfill

\remark \label{_codim_1_CR_Remark_}
Given a real hypersurface $M\subset P$ in a complex
manifold, $\dim_\C P=n$, the rank of $TM\cap I(TM)$ is
$n-1$ everywhere, hence $M$ is a CR-submanifold.

\hfill

Given a CR-manifold $(M,H)$, consider the
bundle $H\oplus \bar H\subset TM\otimes \C$.
This bundle is preserved by the complex conjugation,
hence it is the complexification of a $H_\R\subset TM$.
Since $H\cap \bar H=0$, the map
$\Re:\; H \arrow H_\R$ is an isomorphism.
The map $\1\Id_H:\; H\arrow H$
defines a complex structure operator $I_H$ on $H_\R$,
$I_H^2 = -\Id_{H_\R}$. Clearly, $H$ is the $\1$-eigenspace
of the $I_H$-action on $H_\R\otimes \C$.

\hfill

\remark
We obtain that a CR-structure on a manifold $M$
can be defined as a pair $(H_\R, I_H)$, where
$H_\R\subset TM$  is a subbundle
 in $TM$,
and $I_H\in \End(H_\R)$ is an endomorphism,
$I_H^2 = -\Id_{H_\R}$, such that the $\1$-eigenspace
of $I_H$-action on $H_\R\otimes \C$ satisfies
\[
[H,H]\subset H
\]
This is the definition we shall use.

\hfill

\definition
A \bftext{CR-holomorphic} function on a CR-manifold
$(M,H)$ is a function $f:\; M \arrow \C$ which satisfies
$D_V(f)=0$ for any $V\in\bar H$ ($D_V$ denotes the
derivative).
A {\bf CR-holomorphic} map is a smooth map of CR-manifolds
such that a pullback of CR-holomorphic functions is
CR-holomorphic.

\hfill

\definition
Let $M$ be a smooth manifold, and $R\subset TM$
a subbundle. Consider the commutator $[R, R]\arrow TM$.
This map is not $C^\infty(M)$-linear. However, its
composition with the projection to $TM/R$ is linear.
It is called \bftext{the Frobenius tensor of the
distribution $R\subset TM$}.

\hfill

\definition
A \bftext{contact manifold} is a smooth manifold $M$
equipped with a codimension $1$  subbundle $R\subset TM$
such that the Frobenius tensor $R\times R \arrow TM/R$
is a nowhere degenerate skew-symmetric $TM/R$-valued
form on $R$. In this case, $R$ is called \bftext{the contact
distribution on $M$}.

\hfill

\remark\label{_Reeb_fie_Remark_}
The bundle $TM/R$ is one-dimensional, hence trivial
(if oriented). A trivialization $\eta$ of $TM/R$
defines a 1-form on $TM$, which is called \bftext{a contact
  form of $M$}. Its differential $d\eta$ is
nowhere degenerate on the contact distribution $R$. A choice of a
trivialization $\eta$ also defines a
\bftext{Reeb vector field} $\xi$ by the conditions:
$\xi\rfloor\eta=1$, $\xi\rfloor d\eta=0$.

\hfill

\remark\label{_conical_symple_Remark_}
Let $(M, R)$ be a contact manifold, and $\eta$
a contact form. Using $\eta$, we define a trivialization
of $TM/R$. Then the total space $S$ of positive vectors in
$(TM/R)^*$ is identified with the cone $M\times \R^{>0}$.
We consider the contact structure as a $TM/R$-valued
1-form on $M$. This gives a canonical 1-form on
$\Tot((TM/R)^*)$.

Let $t$ be a unit parameter on $\R^{>0}$, and $t\theta$ the
corresponding 1-form on $S$. It is easy to check
that $d(t\theta)$ is a symplectic form on $S$.
The converse is also true. Starting from a symplectic
form $\omega$ on a cone $M\times \R^{>0}$,
satisfying $\rho(q)^* \omega=q^2\omega$, where
$\rho(q)(m, t) = (m, qt)$ is a dilatation map,
we may reconstruct the contact structure
on $M$ and the contact form. This can be summarized
by saying that a {\bf contact form on $M$ is the same
as a conical symplectic structure on $M\times \R^{>0}$}.
This construction is explained in greater detail
in most textbooks on contact geometry, \emph{e.g.}
\cite{_Arnold:Mechanics_}.

\hfill

\definition
A \bftext{contact CR-manifold} is a CR-manifold
$(M, H_\R, I_H)$, such that the distribution
$H_\R\subset TM$ is contact.

\hfill

\remark
Given a  CR-manifold $(M, H_\R, I_H)$, with $H_\R$
of codimension 1, the Frobenius 2-form
$H_\R\times H_\R \arrow TM/R$ is of type $(1,1)$
with respect to the complex structure on $H_\R$.
Indeed, this form vanishes on $H$ and $\bar H$, because
$[H, H]\subset H\subset H_\R \otimes \C$.

\hfill

\definition
In these assumptions, the $(1,1)$-form
$H_\R\times H_\R \arrow TM/R$ is called \bftext{the Levi form}
of the CR-manifold $(M, H_\R, I_H)$.

\hfill

\definition\label{_CR_stri_pse_Definition_}
A CR-manifold $(M, H_\R, I_H)$ with $H_\R$
of codimension 1 is called \bftext{pseudoconvex}
if the Levi form $\omega$ is positive or negative,
depending on the choice of orientation. If this
form is also sign-definite, then $(M, H_\R, I_H)$
is called {\bf strictly pseudoconvex}, or
{\bf contact pseudoconvex}.

\hfill

\remark\label{_Stein_CR_Remark_}
Let $S\subset P$ be a Stein domain in a complex manifold,
and $\6S$ its boundary. Assume that $\6S$ is smooth;
then $\6S$ inherits a natural CR-structure from $P$.
It is well known that in this case the Levi form on
$P$ is positive, though not always definite (see \cite{grauert}).

\subsection{Automorphisms of CR-manifolds}

\definition
Let $\phi:\; M\arrow M'$ be a smooth map of CR-manifolds
$(M, H, I)$, $(M', H', I')$.
If $\phi$ maps $H$ to $H'$ and commutes
with the complex structure, $\phi$ is called
{\bf CR-holomorphic}. A CR-holomorphic diffeomorphism
is called {\bf a CR-diffeomorphism}.

\hfill

\definition
Let $(M, H, I)$ be a strictly pseudoconvex CR-manifold.
A vector field $V\in TM$ is called {\bf transversal}
if its image in $TM/H$ is nowhere degenerate.
A diffeomorphism flow on $M$ is called {\bf transversal}
if its tangent field is transversal.

\hfill

The following result is proven in \cite{_Schoen:CR_auto_} (see also
\cite{bgs}):

\hfill

\theorem\label{_CR_auto_compact_Theorem_}
Let $M$ be a compact, strictly pseudoconvex CR-manifold
which is not isomorphic to a sphere with a standard
CR-structure, and $G$ the group of CR-automorphisms of $M$.
Then $G$ is a compact Lie group.

\endproof
\hfill

If $M= S^{2n-1}$ is an odd-\-di\-men\-sional
sphere, the group of CR-\-dif\-feo\-mor\-phisms
of $M$ is isomorphic to $SU(n, 1)$ (see \emph{e.g.} \cite{bgs};
an explicit construction is given in Subsection
\ref{_Sasa_on_sphe_Subsection_}).


\section{Vaisman manifolds and Sasakian geometry}


\subsection{Sasakian manifolds}




\definition\label{_Sasakian_Definition_}
A Riemannian manifold $(M,h)$ of odd real dimension is called
\bftext{Sasakian} if the metric cone $\cac(M)=(M\times \RR^{>0},
t^2h+dt^2)$ is
equipped with a dilatation-invariant complex structure,
which makes $\cac(M)$ a K\"ahler manifold (see \cite{bl},
\cite{bg}).

\hfill

\remark
A Sasakian manifold is naturally embedded as a real hypersurface in
its cone,
$M=(M\times \{1\})\subset \cac(M)$. This defines
a CR-structure on $M$ (\ref{_codim_1_CR_Remark_}).

\hfill

\claim
This CR-structure is contact and pseudoconvex
(that is, strictly pseudoconvex).

\hfill

\noindent{\bf Proof.}
The function $\phi(m,t)=t^2$ on $\cac(M)$
defines a K\"ahler potential on $\cac(M)$
(see e.g. \cite{_Verbitsky:LCHK_}).
The level set of a K\"ahler potential
is strictly pseudoconvex, because its Levi form
is equal to $\6\bar\6\phi\restrict{H}$.
\endproof

\hfill

\remark
It is easy to show that
a Sasakian manifold is equipped with a canonical
contact structure. Indeed, a contact form on
$M$ is the same as a conical symplectic form
on $\cac(M)$, as explained in \ref{_conical_symple_Remark_}.
Such a symplectic form is a part
of the K\"ahler structure on $\cac(M)$.


\hfill

\remark
Let $M$ be a Sasakian manifold. On $M\subset \cac(M)$, consider
the vector field
$\xi=I\left(t\frac d {dt}\right)$ ,
where $t\frac d {dt}$ is the dilatation vector field of the
cone $\cac(M)=(M\times \RR^{>0}, t^2h+dt^2)$, and
$I$ the complex structure operator. Then $\xi\rfloor\eta=1$,
$\xi\rfloor d\eta=0$, hence $\xi$ is the Reeb vector field of $M$.

\hfill

The next result is well known, see for example \cite{bg}:

\hfill

\claim
Let $(M, h)$ be a Sasakian manifold.  The Reeb field is unitary and
Killing: its
flow $\rho(t)$ acts on $M$ by isometries. Moreover, it  preserves
the
CR-structure.

\endproof

\hfill

If the orbits of the Reeb flow of a Sasakian manifold are compact,
then the Sasakian structure is called {\bf quasi--regular}. In this
case, if compact, $M$ fibers in circles over a compact K\"ahler
orbifold. The construction can be reversed if one starts with a
compact Hodge orbifold, cf. \cite[Theorem 2.8]{bg00}.

\hfill

\remark\label{_S^1_inva_Sasakian_Remark_}
Any Sasakian metric $h$ is $S^1$-invariant with respect to
some CR-holomorphic $S^1$-action. Indeed, let $G$ be the closure of the
one-\-pa\-ra\-met\-ric group generated
by the Reeb field. In \cite{ov2} (see also \cite{kami})
it is shown that this group is a
compact torus. Clearly, $h$ is
$G$-invariant. Taking $S^1\subset G$
generated by a vector field
sufficiently close to the Reeb field,
we may also assume that this $S^1$-action
is proper and transversal.

\subsection{Vaisman manifolds}

Our method will be to constantly translate the Sasakian geometry
into locally conformally K\"ahler and Vaisman geometry. Here we recall
the basics.  For details and examples we refer to \cite{drag},
\cite{ov1}, \cite{ov2}, \cite{ov3}, \cite{_Verbitsky:LCHK_}.

As we shall only deal with compact manifolds, we can take as
definition the following characterization:

\hfill

\definition\label{def_vai} A compact complex manifold $(N,I)$ is
called a \bftext{ Vaisman manifold}
if it admits a K\"ahler covering $\Gamma\rightarrow(\tilde
N, I, h)\rightarrow (N,I)$ such that:
\begin{itemize}
\item $\Gamma$ acts by holomorphic homotheties with respect to $h$
(this says that $(N,I)$ is equipped with a {\bf locally conformally
K\"ahler structure}).
\item $(\tilde N, I, h)$ is isomorphic to a K\"ahler cone over a
compact Sasakian manifold $M$. Moreover, there exists a Sasakian
automorphism $\phi$ and a positive number $q>1$ such that $\Gamma$
is isomorphic to the cyclic group generated by $(x,t)\mapsto
(\phi(s), tq)$.
\end{itemize}

\hfill

In particular, the product of a compact Sasakian manifold with $S^1$
is equipped with a natural Vaisman structure.

The K\"ahler metric $h$ on $\cac(M)=M\times \R^{>0}$
has a global potential $\psi$, which is
expressed as $\psi(m,t) =t^2$. The
metric $\psi^{-1}\cdot h$ projects on $N$ into a Hermitian,
locally conformally K\"ahler metric, say $g$, whose
fundamental two-form  $\omega$
satisfies the equation $d\omega=\theta\wedge \omega$ for a closed
one-form $\theta$ called the Lee form. Then $\psi=\mid\theta\mid^2$.

Further on, $N$ is considered as a Hermitian manifold,
with the Hermitian metric $g$. Let $\theta^\sharp$ be the vector
field on $N$ dual to the Lee form $\theta$. Then
$\theta^\sharp$ is called {\bf the Lee field on $N$}.

\hfill

The Lee field $\theta^\sharp$ is Killing, parallel
with respect to the Levi-Civita connection on $N$ and holomorphic.
It thus determines two foliations on $N$:
\begin{itemize}
\item ${\cal F}_1$, one-dimensional, tangent to $\theta^\sharp$.
\item ${\cal F}_2$, holomorphic two-dimensional, tangent to
$\theta^\sharp$ and $I\theta^\sharp$.
\end{itemize}

\remark\label{autom_pot} As $\theta^\sharp$ is parallel, $\mathrm{Lie}_{\theta^\sharp}g(\theta^\sharp,\theta^\sharp)=0$, hence the flow
of $\theta^\sharp$ preserves the potential $\psi$.

\hfill

\proposition If the foliation ${\cal F}_1$ (resp. ${\cal F}_2$) is
quasi--regular (thus having compact leaves), then the leaf space
$N/{\cal F}_1$ (resp. $N/{\cal F}_2$  is a Sasakian (resp.
projective K\"ahler) orbifold.

\hfill

\remark\label{weight_bundle} Let $L$ be the weight line bundle
associated to the Vaisman manifold \emph{via} the subjacent
l.c.K. structure. The Lee form can be interpreted as a canonical
Hermitian connection in its complexification (that we also denote by
$L$) and one can prove (see \cite{_Verbitsky:LCHK_}) that the
curvature is positive except on the Lee field. The Chern connection
in $L$ is trivial along ${\cal F}_2$. Therefore, if $N$ is
quasi-regular, $L$ is a pullback of a Hermitian
line bundle $\pi_* L$ on the K\"ahler orbifold
$N/{\cal F}_2$. Since the projection
$\pi:\; N \arrow N/{\cal F}_2$ kills the directions on which the
curvature was non-positive, the push-forward bundle $\pi_* L$
is ample (cf. \cite{ov2}).

\hfill

The following result from \cite{ov2} will be important in the
sequel:

\hfill

\theorem\label{def} \cite[Proposition 4.6]{ov2}
A compact Vaisman manifold can be deformed into
a quasi--regular Vaisman manifold, with the same K\"ahler covering
$(\tilde N, h)$.

\endproof

\hfill


%


\section{Algebraic cones and CR geometry}


\subsection{Algebraic cones}
\label{_cones_Section_}

\definition
Let $X$ be a projective
variety, and $L$ an ample line bundle on $X$. The \bftext{algebraic
cone} $\cac(X,L)$ of $X$ is the total space
of non-zero vectors in $L^*$.
A \bftext{cone structure} on $\cac(X,L)$ is the $\C^*$-action
arising this way (by fibrewise multiplication).

\hfill

\definition\label{_alge_cone_clo_ope_Definition_}
Let $\cac(X, L)$ be an algebraic cone. Consider
the associated affine variety $\bar\cac(X,L):=\Spec \oplus_i H^0(X,
L^i)$.
Geometrically, $\bar\cac(X,L)$ is a complex
variety, obtained by adding a point at the ``origin''
of the cone $\cac(X, L)$. We call $\bar\cac(X,L)$
\bftext{the closure} of the algebraic cone $\cac(X, L)$.
This space is called {\bf a closed algebraic cone},
and $\cac(X, L)$ {\bf an open algebraic cone}.

\hfill

The definition of an algebraic cone
is motivated by the following observation.
Let $h$ be a Hermitian metric on $L^*$, such that
the curvature of the associated Hermitian connection
is negative definite (such a metric exists, because
$L$ is ample). Consider a function
$\cac(X,L)\stackrel\phi\arrow \R$, $\phi(v)=h(v,v)$.
Then $\6\bar\6\phi$ is a K\"ahler metric on $\cac(X,L)$
(see e.g. \cite{_Verbitsky:LCHK_}).
The associated K\"ahler manifold is a Riemannian cone
of a unit circle bundle
\[
\{ v\in \cac(X,L)\ \ |\ \ h(v,v)=1\}
\]
which is, therefore, Sasakian.

\hfill

In the Introduction, we defined algebraic cones
in a less constructive manner (see \ref{_algebra_cone_intro_Definition_}).
This definition is equivalent to the one given above,
as follows from \cite{ov2} and \cite{ov3}.
Starting from an algebraic cone in the sense
of \ref{_algebra_cone_intro_Definition_}, that is,
an affine algebraic variety $X$ with an action $\rho$ of
$\C^*$ contracting $X$ to a single singular point $x_0$,
we may embed $(X\backslash x_0)/\rho(2)$
into a diagonal Hopf manifold, as shown in
\cite{ov3}. This allows us to equip
$(X\backslash x_0)/\rho(2)$ with a Vaisman
metric. In \cite{ov2},
it was shown that a covering of a compact Vaisman
manifold is isomorphic to a space of
non-zero vectors in some anti-ample line bundle
over a projective orbifold. Therefore, it is an open algebraic cone,
in the sense of \ref{_alge_cone_clo_ope_Definition_}.
This implies that \ref{_algebra_cone_intro_Definition_}
is equivalent to  \ref{_alge_cone_clo_ope_Definition_}.

\hfill

The arguments of the present paper are built on the
correspondence between the Sasakian manifolds and the
algebraic cones, which is implied by the following
proposition.

\hfill

\proposition
Let $M$ be a compact Sasakian manifold,
$\cac(M)$ its cone, considered as a complex manifold.
Then $\cac(M)=\cac(X,L)$ is an algebraic
 cone, associated with a projective orbifold $X$.

\hfill

\noindent{\bf Proof.}
Indeed, the product $M\times S^1$ is Vaisman
(see above), and is covered by the K\"ahler
cone $\cac(M)$.
Proposition 4.6 of \cite{ov2} implies that the same cone
$\cac(M)$ is a covering of a quasi-regular Vaisman manifold,
that is, a total space of an elliptic fibration
$E\arrow X$, with $X$ a projective orbifold.
But any quasi-regular
Vaisman manifold can be obtained as a quotient
of a cone $\cac(X,L)$ (see \ref{weight_bundle}) by an equivalence
$t\sim qt$, where $q\in \C^*$ is a fixed
complex number, $|q|>1$ (see \ref{def_vai}).
 \endproof

\hfill

The K\"ahler structure of the Riemannian cone $\cac(M)$ explicitly
depends on  the Sasakian metric of $M$. However, the holomorphic
structure of the cone is determined by the underlying CR-structure of
$M$, as follows from \ref{_Sasa_unique_Theorem_}.
A weaker version of this statement is obtained immediately
from standard results of complex analysis.

\hfill

\proposition\label{_cac_lambda_unique_Proposition_}
Let $M$ be a compact Sasakian manifold,
and $\lambda$ a positive real number. Denote by $\cac(M)_\lambda$
the
set of all $(m,t)\in \cac(M)$, $t\leq \lambda$.
Then the holomorphic structure of $\cac(M)_\lambda$ depends only
on the CR-structure of $M$.

\hfill

\noindent{\bf Proof.} Consider the standard embedding
$M\hookrightarrow \cac(M)$, $m \arrow m\times \{\lambda\}$,
and let $M_\lambda$ be its image, which is the boundary
of a complete Stein domain $M\times [0,\lambda]/M\times\{0\}$
(see \cite{ov3}), Theorem 3.1).
Let $V_\lambda$ be the space of CR-holomorphic functions on
$M_\lambda$.
Using the solution of  $\bar\6$-Neumann problem
(\cite{_Marinesku_Dinh_}), we
identify $V_\lambda$ with the space of holomorphic
functions on the Stein domain $\cac(M)_\lambda$
which are smooth on its boundary $M_\lambda$. Then
$\cac(M)_\lambda$ is the holomorphic spectrum
of the ring $V_\lambda$.
\endproof

\hfill

The following problem is then natural: if one starts with an
algebraic cone and fixes a pseudoconvex CR-hypersurface in it, when
does the CR-structure underlie a Sasakian structure? As
shown in the next section, the answer is related to the K\"ahler
potentials on the cone.

\subsection{Sasakian manifolds in algebraic cones}

\theorem\label{_Sasa_in_cone_Theorem_}
Let $M$ be a smooth real hypersurface in a closed
algebraic cone $\cac$,
considered as a CR-manifold. Assume that $M$ is contact
and pseudoconvex (this implies that $M$ is the boundary
of a Stein domain $\cac_1$ in  $\cac$).
Then $M$ admits a Sasakian metric if and only if
for some cone structure $\rho:\; \C^*\arrow \Aut(\cac)$,
$M$ is $S^1$-invariant.

\hfill

\noindent{\bf Proof.} The "if" part follows
{}from \ref{_S^1_inva_Sasakian_Remark_}, where
the $S^1$-action is constructed.
For the converse, assume that $\cac_1$ is an open,
$S^1$-invariant subset of a cone.
To prove that its boundary
$M$ is Sasakian, we need to construct
a K\"ahler potential which is homogeneous under
the cone action and such that $M$ is its level set. To this end, we
introduce the following

\hfill

\definition \bftext{A section} of the action
$\rho:\; \R^{>0} \arrow \Aut(\cac)$ is a subset $V\subset \cac$
such that $\rho(\lambda_1) V$ does not intersect
$\rho(\lambda_2) V$ for $\lambda_1\neq \lambda_2$,
and $\rho(\R^{>0})V=\cac$.

\hfill

We now fix a K\"ahler potential $\psi$ on $\cac_1$
mapping its boundary to 1 (the existence of such a potential is
assured by the strict pseudoconvexity of $M$; see \emph{e.g.}
\cite[Consequence 3.2]{_Marinesku_Yeganefar_}).
Let $\Delta_1$ be a unit disk in $\C$. Averaging
$\psi$ with $S^1$-action induced by $\rho$, we may assume that
$\psi$ is $S^1$-invariant.  For all $m\in M$, the
discs $\rho(\Delta_1)m$, bounded by the images of $S^1$, belong
to $\cac_1$ (indeed, being strictly
plurisubharmonic, $\psi$ is subharmonic
on all curves in $\cac$).
This implies that $M\subset \cac$ is a section of
$\rho:\; \R^{>0} \arrow \Aut(\cac)$, in the sense of the above
definition.

This allows us to define a map
$\phi:\; \cac \arrow \R^{>0}$,
mapping $x\in \cac$ to $t^2$, where
$x \in \rho(t)M$. By construction, this map
is homogeneous with respect
to the cone action.

To finish the proof it remains to show that
$\phi$ is a K\"ahler potential.
Clearly, on the contact distribution
of $M$ the form $\6\bar\6 \phi$ is equal
to $\phi \omega_0$, where $\omega_0$
is the Levi form, and it is
positive because $M$ is pseudoconvex.
On the plane generated by $\rho$-action,
$\phi=|z|^2$, hence  plurisubharmonic.
Finally, these two spaces are orthogonal,
because $\phi$ is $\rho(S^1)$-invariant.
Then $\phi$ is plurisubharmonic. \endproof

\hfill

\remark
In the proof of
\ref{_Sasa_in_cone_Theorem_},
we constructed an {\em $S^1$-invariant}
Sasakian metric on $M$. Moreover, the Reeb
field of $M$ is proportional (with constant
coefficient) to the tangent field to the $S^1$-action.
However, $M$ can
possibly admit other Sasakian metrics, not all of them
necessarily $S^1$-invariant or having the prescribed
Reeb field; see \cite{bgs} or \ref{_when_CR_admits_Sasakian_Theorem_}.


\section{Existence and uniqueness of Sasakian structures}
\label{_Sasa_exi_Section_}


\subsection{Sasakian structures on CR-manifolds with $S^1$-action}\label{exist_subsection}

\theorem\label{_Sasa_if_S^1_acts_Theorem_}
 Let $M$
be a  compact pseudoconvex contact
CR-manifold, $\dim M\geq 5$. Then the following assumptions are equivalent.

\begin{description}
\item[(i)] $M$ admits a transversal,
CR-holomorphic action of $S^1$.

\item[(ii)]  $M$ admits a transversal,
CR-holomorphic vector field.

\item[(iii)] $M$ admits a Sasakian metric, compatible with
the CR-structure.
\end{description}

\hfill

\noindent{\bf Proof:} The implication (i) $\Rightarrow$ (ii)
is clear, and (iii) $\Rightarrow$ (i) follows from
\ref{_S^1_inva_Sasakian_Remark_}. The implication
(ii) $\Rightarrow$ (i) is clear from \ref{_CR_auto_compact_Theorem_}.
Indeed, since the group of CR-automorphisms of $M$ is
compact (unless $M$ is a sphere), any diffeomorphism
flow can be approximated by an $S^1$-action within
its closure. The implication (ii) $\Rightarrow$ (iii)
follows from differential-geometric arguments (see
\ref{_another_arguments_citation_Remark_}). We use
another argument, which is based on complex analysis.
The implication (i) $\Rightarrow$ (iii), and
the proof of \ref{_Sasa_if_S^1_acts_Theorem_},
is given by the following proposition.

\hfill

\proposition \label{_embe_to_cone_CR_Sasakian_Proposition_}
Let $M$ be a compact pseudoconvex contact
CR-manifold, $\dim M\geq 5$. Assume that $M$ admits
a proper CR-holomorphic $S^1$-action $\rho$. Then
$M$ admits an $S^1$-equivariant CR-embedding to an
algebraic cone.

\hfill

\noindent{\bf Proof:}
As shown in \cite{_Marinesku_Dinh_}, $M$
is the boundary of a Stein variety $\cac_1$ with isolated
singularities.
Since $\cac_1$ is constructed uniquely
(by solving the boundary $\bar\6$-Neumann problem),
the $S^1$-action on $M$ can be extended to a
holomorphic $S^1$-action on $\cac_1$. We shall
explain now how to integrate $\rho$ to a
$\C^*$-action.

Let $\zeta$ be the tangent vector field induced by the
$S^1$-action, and $\zeta^c= I(\zeta)$
the corresponding vector field in $T\cac_1$.
Since $\zeta$ is transversal with respect to CR-structure,
$\zeta^c$ has to point inward or outward,
everywhere in $M$ (in other words, it has to point either
towards the filled-in part or towards the
opposite direction). Assume it is inward.

\begin{wrapfigure}[16]{r}{13em}
\epsfig{file=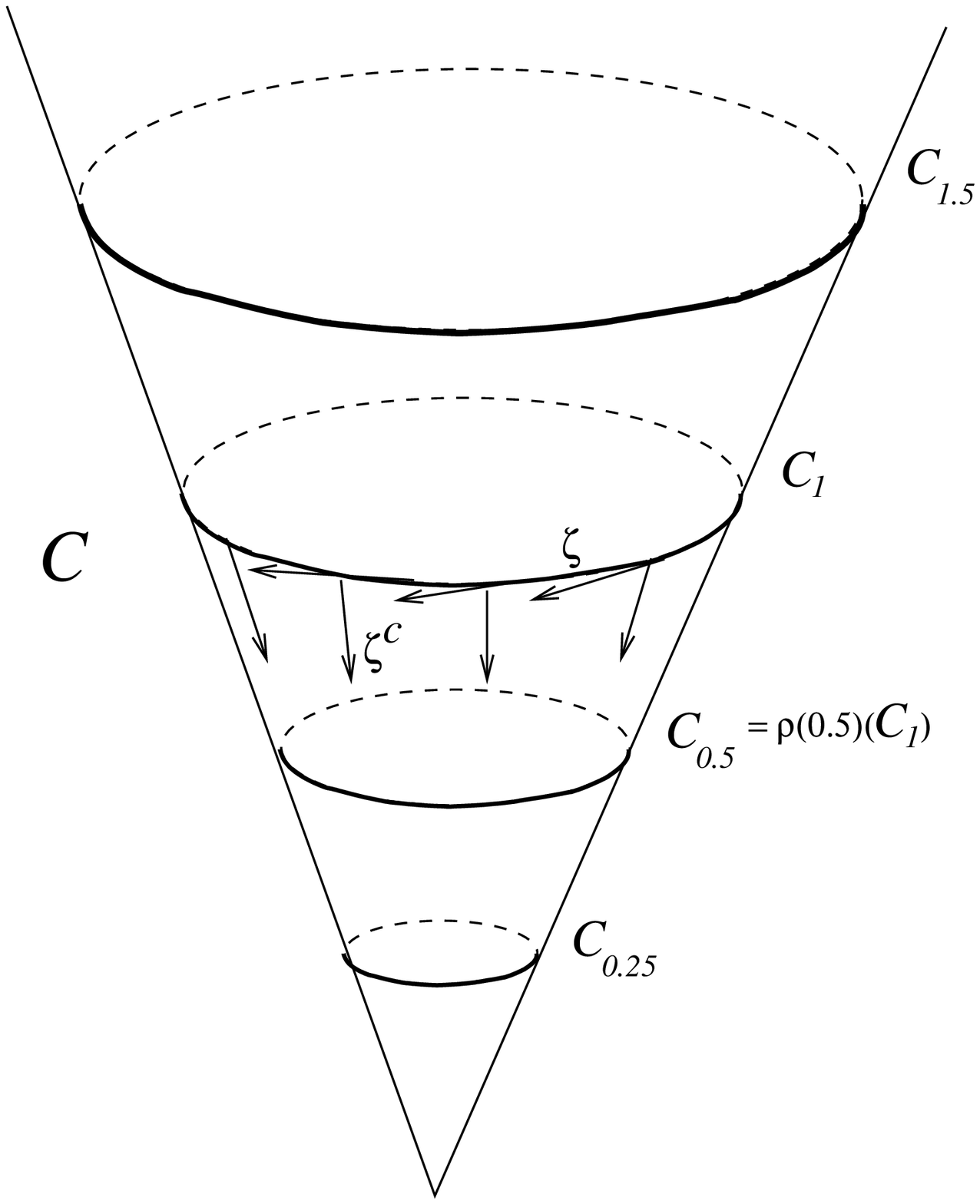,width=13em}
\caption{Gluing $\cac_1$ to itself}
\end{wrapfigure}

The flow of  $\zeta^c$ provides
a holomorphic automorphism $\rho_1$ of $\cac_1$,
mapping $\cac_1$ into itself.
Iterating this map, we find that we can integrate
$\rho$ from $S^1$ to an action of $0< |z| \leq 1$.
Inverting this construction and gluing images
of $\rho(0< |z| <\epsilon)$ together, as shown in
Figure 1,  we obtain a domain $\cac$ containing $\cac_1$
where $\rho$ can be integrated to a
$\C^*$-action.

The potential $\phi$, constructed in
the proof of  \ref{_Sasa_in_cone_Theorem_}, can be defined
in the same way now, and it gives a homogeneous, $S^1$-invariant
plurisubharmonic function on $\cac$ as indicated.
Then $\cac/\rho(2)$ is Vaisman, hence $\cac$ is
a cone of a projective orbifold
as follows from \cite[Proposition 4.6]{ov2}. \endproof

\hfill

This finishes the proof of \ref{_S^1_equiv_Sasakian_Theorem_}.
The proof of \ref{_embe_to_cone_CR_Sasakian_Proposition_}
also leads to the following result.

\hfill

\corollary \label{_embe_to_cone_CR_Unique_Corollary_}
Let $M$ be a compact pseudoconvex contact
CR-manifold, admitting a proper CR-holomorphic
$S^1$-\-ac\-tion $\rho$, and \[ M\hookrightarrow \cac\]
an $S^1$-equivariant embedding to an algebraic
cone. Then $\cac$ is uniquely determined
by $M$ and the $S^1$-action.

\hfill

\noindent{\bf Proof:} By \ref{_cac_lambda_unique_Proposition_},
the holomorphic structure on  $\cac_1$ is determined
by the CR-structure on $M$. The cone $\cac$ is reconstructed
from $\cac_1$ and the $S^1$-action as above.
\endproof

\subsection{Uniqueness of the algebraic cone}\label{uni_subsection}

Every Sasakian manifold
is CR-embedded into an algebraic cone by
\ref{_embe_to_cone_CR_Sasakian_Proposition_}.
This cone is determined uniquely by an $S^1$-action,
as follows from \ref{_embe_to_cone_CR_Unique_Corollary_}.
On the other hand, the cone is determined
by the Sasakian metric. Then, if
a given Sasakian metric is invariant under two different
$S^1$-actions, the cone associated to one $S^1$-action
is isomorphic to the cone associated to the other $S^1$-action.
Therefore, unless $M$ is a sphere, \ref{_Sasa_unique_Theorem_}
is implied by the following proposition.

\hfill

\proposition\label{_CR_auto_inva_Sasa_Proposition_}
(\cite[Proposition 4.4]{bgs}) Let $M$ be a CR-manifold of Sasakian type,
and $G$ the group of CR-automorphisms of $M$. Assume that $M$ is not
a sphere. Then $M$ admits a $G$-invariant Sasakian metric $g$.

\hfill

\endproof


\section{The positive Sasakian cone of a CR-manifold}


\subsection{Positive Sasakian cone}
\label{_posi_Sasa_cone_Subsection_}

The notion of  positive Sasakian cone is due to \cite{bgs}.
We use it to classify the Sasakian metrics on a sphere.

\hfill

\definition\label{_posi_sasakian_cone_Definition_}
Let $(M, H)$ be a strictly pseudoconvex CR-manifold. The
Levi form $H\otimes H \arrow TM/H$ is sign-definite.
This gives an orientation on $TM/H$. A transversal CR-holomorphic
vector field is called
\bftext{positive} if its projection to $TM/H$
is everywhere positive. The \bftext{positive Sasakian cone}
is the space of all transversal, positive, CR-holomorphic
vector fields.

\hfill

\definition
Let $M$ be a Sasakian manifold, $\cac(M)= M\times {\Bbb R}^{>0}$
its cone, with $t$ a coordinate in $\Bbb R{>0}$,
and $\frac{d}{dt}$ the corresponding holomorphic
vector field. It is clear from the definition
that $\xi:= I(t\frac{d}{dt})$ is tangent to the fibration
$M\times \{t\}$, hence defines a vector field
on $M$. We normalize it in such a way that
$|\xi|=1$. Then $\xi$ is called \bftext{the Reeb
field of the Sasakian manifold $M$}.

\hfill

This definition is compatible with the one
used in contact geometry (\ref{_Reeb_fie_Remark_}).
The following well-known claim is easy to prove (see, for example, \cite{bg}):

\hfill

\claim
In these assumptions, $\xi$ is transversal,
positive, CR-holomorphic and Killing.

\hfill

The following theorem is implied by
Lemma 6.4 of \cite{bgs}.

\hfill

\theorem\label{_Sasaki_cone_Theorem_}
Let $M$ be a compact, strictly pseudoconvex
CR-manifold, and ${\cal R}$ the Sasakian
positive cone. Denote by ${\cal S}$ the set
of Sasakian metrics on $M$, and let
${\cal S} \stackrel \Psi \arrow {\cal R}$ map a metric
into the corresponding Reeb field.
Then $\Psi$ is a bijection.

\hfill

\noindent{\bf Proof:} We give an independent proof
of \ref{_Sasaki_cone_Theorem_}, using the
same kind of arguments as we used in the proof
of \ref{_S^1_equiv_Sasakian_Theorem_}.
The Reeb field on a contact, pseudoconvex
CR-manifold $(M, H, I)$ determines the Sasakian
metric uniquely, as can be seen from the
following argument. Denote the corresponding
contact form by $\eta$ (see \ref{_Reeb_fie_Remark_}).
The Hermitian form $d\eta\restrict H$
is equal to the Levi form by construction;
the Reeb field $\xi$ is orthogonal to $H$
and has length 1. Therefore, the map
 ${\cal S} \stackrel \Psi \arrow {\cal R}$
is injective. It remains to show that $\Psi$
is a surjection.

Let $\zeta$ be a positive, transversal
CR-holomorphic vector field on a compact, strictly pseudoconvex
CR-manifold $M$, and $B$ the corresponding Stein domain,
$\6 B = M$. Then $e^{t\zeta}$ induces an automorphism
of $B=\Spec(\calo_M)$, where $\calo_M$ is the ring
of CR-holomorphic functions on $M$. Since $\zeta$
is positive, the vector field $-\zeta^c:= - I(\zeta)$
points transversally to $M$ towards $B$. Therefore, the map
$e^{-t\zeta^c}:\; B \arrow B$ is well defined for small $t$, and
maps $B$ to a strictly smaller subset which
is contained in the interior of $B$. Iterating
this map, we obtain that $e^{-t\zeta^c}$ is
well defined for all $t$. Inverting this
procedure as in the proof of \ref{_embe_to_cone_CR_Sasakian_Proposition_},
we obtain that $\zeta^c$ induces a holomorphic action $\rho$ of
the multiplicative group $\R^{>0}$
on a Stein domain $B_\infty$, which contains $M$
as a hypersurface. Clearly, $\rho$ is a contraction,
with $\rho(\epsilon_i)$, $\epsilon_i \arrow 0$,
putting $B$ into a sequence of smaller open
balls converging to a single fixed point $x_0$
(see \cite{ov3}). Therefore, $\rho$
is free outside of $\{ x_0\}$, and
\[
  B_\infty \backslash S= \rho(\R^{>0}) M.
\]
because each orbit of $\rho$ encounters $M$ on the way to $x_0$.
We define a function $\phi:\; (B_\infty \backslash \{x_0\}) \arrow \R^{>0}$
by
\[
\phi(x)= \lambda^2, \ \ \text{for}\ \ \rho(\lambda^{-1})x \in M.
\]
Then $\phi$ is a K\"ahler potential on $B_\infty$. Indeed,
on the contact distribution $H\subset TM$, $\6\bar\6\phi$
is proportional to the Levi form of $M$, because
$\phi$ is constant on $M$. In particular,
$\6\bar\6\phi$ is positive on $H$.
On the 1-dimensional complex foliation $F$ generated
by $\langle\zeta, \zeta^c\rangle$, $\phi$ is quadratic,
and can be written in appropriate holomorphic
coordinates as $z\arrow |z-c|^2$. Finally,
 $\6\bar\6\phi$ vanishes on pairs $(x, y)$,
$x\in H, y\in F$, because $\Lie_{\zeta}\phi=0$.
Therefore, $\6\bar\6\phi$ is positive
on $F\oplus H= T B_\infty$.
The function $\phi$ is by construction homogeneous with respect
to the action of $\rho$. Therefore, the corresponding
K\"ahler form $\6\bar\6\phi$ is also homogeneous,
and $B_\infty$, considered as a Riemannian
manifold, is identified with a Riemannian cone over $M$.
This gives a Sasakian metric on $M$. It is easy
to check that the corresponding Reeb field
is equal to $\zeta$. We proved that any
positive transversal CR-holomorphic vector field
is induced by some Sasakian metric.
\ref{_Sasaki_cone_Theorem_} is proven.
\endproof

\hfill

Now we can prove \ref{_when_CR_admits_Sasakian_Theorem_}.
A Sasakian manifold $M$ is \bftext{quasi-regular}
if the 1-dimensional foliation $F_1$ induced by the
Reeb field on $M$ has compact fibers. Quasiregular
Sasakian manifolds are always obtained from the
construction described in \ref{_U(1)_fibration_Examples_}
(see \cite{bg00}). Therefore, to prove \ref{_when_CR_admits_Sasakian_Theorem_},
we need to show that a given CR-manifold $M$ admits
a quasi-regular Sasakian structure, if it admits
some Sasakian structure.

Denote by $A_0$ the 1-parameter group
of CR-holomorphic isometries, generated by $e^{t\xi}$,
where $\xi$ is the Reeb field of $M$.
Let $A$ be its closure in the Lie
group of CR-holomorphic isometries of $M$.
Since $A_0$ is abelian, $A$ is
also abelian; it is compact, because
the group of isometries is compact. Therefore,
$A$ is a compact torus.

Let ${\goth a}$ be its Lie algebra.
We consider ${\goth a}$ as a subset in
the space of CR-holomorphic vector fields
on $M$. We call a vector field $\zeta\in TM$
\bftext{quasi-regular} if  the corresponding 1-dimensional
foliation has compact fibers. A vector field
$\zeta\in {\goth a}$ is quasi-regular if
it is tangent to an embedding $S^1\hookrightarrow A$.
Such embeddings correspond to rational
points in ${\goth a}$, hence they are dense
in ${\goth a}$. Taking a quasiregular
$\zeta\in {\goth a}$ sufficiently close to
$\xi$, we may assume that it is
also transversal and positive.
By \ref{_Sasaki_cone_Theorem_},
the corresponding Sasakian manifold
is quasiregular. We proved
\ref{_when_CR_admits_Sasakian_Theorem_}.

\hfill

\remark\
A similar deformation-type argument was used in
 \cite[Proposition 1.10]{kami}.

\subsection{Sasakian metrics on a sphere}
\label{_Sasa_on_sphe_Subsection_}

To finish the proof
of \ref{_Sasa_unique_Theorem_},
we need to consider the case of a sphere.
Let $M= S^{2n-1}\subset \C^n$ be an odd-dimensional
sphere equipped with a standard CR-structure.
We are going to classify the Sasakian metrics
compatible with this CR-structure. We are interested
in Sasakian metrics up to CR-automorphism.

Let $G$ be a group of CR-automorphisms of $M$.
It is well known that $G\cong SU(n,1)$
(see e.g. \cite{bgs}). Using the same argument
as used in the proof of \ref{_cac_lambda_unique_Proposition_},
we may assume that $G$ acts as a group of holomorphic
automorphisms on an open ball $B\subset \C^n$, $\6 B = M$.
The action of $SU(n,1)$ on $B$ is very easy to describe
explicitly. Let us identify $B$ with a projectivization
of the positive cone
\[
\{ \xi \in V \ \ |\ \ (\xi, \xi)_p >1\},
\]
where $(\cdot, \cdot)_p$ is a Hermitian form of
signature $(n, 1)$ on $V=\C^{n+1}$. The group $U(n, 1)$ acts
on $V$ preserving the metric, hence
$SU(n,1)\subset PU(n, 1)$ acts on $B\subset {\Bbb P}V$.
This action is holomorphic, therefore its restriction to
$M= \6 B$ is CR-holomorphic. Clearly, $G= SU(n, 1)$
acts on the interior of $B$ transitively. This gives

\hfill

\proposition\label{_CR_transi_on_ball_Proposition_}
Let $M\subset \C^n$ be an odd-dimensional sphere,
considered as a CR-manifold, $G= \Aut_{CR}(M)$
the group of CR-automorphisms, $G\cong SU(n,1)$.
Then $G$ acts transitively on the interior part of the
open ball $B$, $M= \6 B$.

\endproof

\hfill

We are interested in classification of Sasakian structures
up to CR-automorphism. A Sasakian structure
on $M$ induces a $\C^*$-action on a Stein domain
containing $B$ (\ref{_Sasa_unique_Theorem_}).
This way, $B$ is identified with an open part of an
algebraic cone. Since $G$ acts on $B$ transitively,
it maps the origin of this cone into any other
interior point of $B$.

\hfill

\corollary\label{_may_fix_0_Sasakian_on_sphere_Corollary_}
In assumptions of \ref{_CR_transi_on_ball_Proposition_},
let $g$ be a Sasakian metric on $M$, and $\xi$ its
Reeb field. As shown in Subsection
\ref{_posi_Sasa_cone_Subsection_},  $\rho(t):=e^{-t\xi^c}$
acts on $B$ by holomorphic contractions. Denote by
$x_0$ the fixed point of $\rho$.
Then, after an appropriate action of
the group of CR-automorphisms of $G$, we
may assume that $x_0$ is $0\in B$.

\endproof

\hfill

The group $G_0$ of CR-automorphisms of $M$ fixing $0\in B$
is identified with the stabilizer of $0$ under
$SU(1,n)$-action on $B$, that is, with $U(n)$.
Denote by ${\cal R}_0$ the part of positive
Sasaki cone consisting of those positive
transversal CR-holomorphic vector fields
$\xi$ which fix $0\in B$. As follows
from \ref{_Sasaki_cone_Theorem_} and
\ref{_may_fix_0_Sasakian_on_sphere_Corollary_},
every Sasakian metric on $M$ corresponds
to some $\xi\in {\cal R}_0$, up
to a CR-automorphism. Then, the set ${\cal S}/G$
of isomorphism classes of Sasakian metrics
is identified with
${\cal R}_0/G_0$

Clearly,
${\cal R}_0$ is the set of all
$\xi\in {\goth u}(n)$ which
are positive and transversal, that is,
have all eigenvalues $\alpha_i$
with $-\1\alpha_i$ positive real numbers.
The group $U(n)$ acts on ${\cal R}_0$
in a natural way, and each orbit
is determined by the corresponding set
of eigenvalues. This gives the following theorem:

\hfill

\theorem\label{_Sasakian_on_sphere_classifi_Theorem_}
Let $M\subset \C^n$ be an odd-dimensional sphere,
considered as a CR-manifold, $G= \Aut_{CR}(M)$
the group of CR-automorphisms, $G\cong SU(n,1)$,
and ${\cal S}$ the set of Sasakian metrics on
$M$. Then ${\cal S}/G$ is in natural,
bijective and continuous correspondence with
the set of unordered $n$-tuples of positive
real numbers.

\endproof

\hfill

{}From this construction, it is clear that the Riemannian
cone of each Sasakian structure on a sphere is identified
naturally with $\C^n\setminus\{0\}$. This proves \ref{_Sasa_unique_Theorem_}
in the case when $M$ is a sphere. We finished the proof
of \ref{_Sasa_unique_Theorem_}.

\hfill

\hfill

\noindent{\bf Acknowledgements:} The authors thank Charles
Boyer for an illuminating discussion on a first draft of the paper and
Sorin Dragomir for useful informations about Lee's and Webster's work in
pseudo-Hermitian geometry. The authors thank Max Planck Institut in Bonn for hospitality
during part of the preparation of this paper.

{\small

\noindent {\sc Liviu Ornea\\
University of Bucharest, Faculty of Mathematics, \\14
Academiei str., 70109 Bucharest, Romania.}\\
\tt Liviu.Ornea@imar.ro, \ \ lornea@gta.math.unibuc.ro

\hfill

\noindent {\sc Misha Verbitsky\\
University of Glasgow, Department of Mathematics, \\15
  University Gardens, Glasgow, Scotland.}\\
{\sc  Institute of Theoretical and
Experimental Physics \\
B. Cheremushkinskaya, 25, Moscow, 117259, Russia }\\
\tt verbit@maths.gla.ac.uk, \ \  verbit@mccme.ru
}

\end{document}